\documentclass[reqno,12pt]{amsart}

\usepackage{epsf}
\usepackage{graphics}
\usepackage{amssymb}
\usepackage{amsmath}

\date{}

\theoremstyle{plain}
\newtheorem{theorem}{Theorem}
\newtheorem{corollary}{Corollary}
\newtheorem{lemma}{Lemma}
\newtheorem{proposition}{Proposition}

\theoremstyle{definition}

\theoremstyle{remark}

\newtheorem*{remark}{Remark}

\def\Z{{\mathbb Z}}

\def\R{{\mathbb R}}

\title{Trefoil Plumbing} 

\author{Sebastian Baader \and Pierre Dehornoy}

\begin{document}

\begin{abstract} We give a criterion for an open book to contain an $n$-times iterated Hopf plumbing summand. As an application, we show that fibre surfaces of positive braid knots admit a trefoil plumbing structure.
\end{abstract}

\maketitle

\section{Introduction}

An open book is a connected orientable surface $\Sigma$ with boundary, together with a diffeomorphism $\phi: \Sigma \to \Sigma$ fixing the boundary pointwise. Positive Hopf plumbing is an operation that adds a $1$-handle $h$ to $\Sigma$ and composes $\phi$ with a right-handed Dehn twist along an embedded circle that runs once through $h$. In case an open book $(\Sigma,\phi)$ represents the $3$-sphere, the boundary of $\Sigma$ is a fibred link with monodromy~$\phi$. Positive Hopf plumbing was originally defined in this classical setting by Stallings~\cite{St}, where it corresponds to gluing a positive Hopf band on top of a fibre surface.
In this paper, we consider positive $n$-plumbing, an $n$-times iterated version of positive Hopf plumbing, as sketched in Figure~1, for $n=4$.\footnote{This is in fact a schematic picture indicating the positions of $n$ successive $1$-handles. In order to get an embedded picture, one would have to add a positive full twist to every $1$-handle.}
\begin{figure}[ht]
\scalebox{1.0}{\raisebox{-0pt}{$\vcenter{\hbox{\epsffile{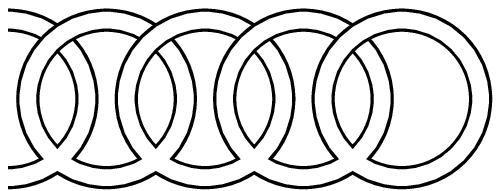}}}$}} 
\caption{}
\end{figure}

The special case $n=2$ is called positive trefoil plumbing, since it amounts to plumbing the fibre surface of a positive trefoil knot on top of a surface. We say that an open book admits a positive $n$-plumbing summand, if it is obtained from a suitable open book by one positive $n$-plumbing operation. Open books that arise from the trivial open book $(D^2,\text{Id})$ by iterated positive Hopf plumbing are known to be right-veering~\cite{HKM}. Our first result is a plumbing criterion for such open books.

\begin{theorem} An open book $(\Sigma,\phi)$ with right-veering monodromy $\phi$ admits a positive $n$-plumbing summand, if and only if there exists a relatively embedded arc $I \subset \Sigma$, such that $I,\phi(I),\phi^2(I),\ldots,\phi^n(I)$ are homologically independent and pairwise disjoint, except at their common boundary points, after a suitable relative isotopy.
\end{theorem}

This is a generalisation of the Hopf deplumbing criterion proved by Etnyre and Li~\cite{EL} (Theorem~3.3). In terms of the homological arc complex defined in their paper, the last condition means that the arcs $I,\phi(I),\phi^2(I),\ldots,\phi^n(I)$ form an ($n+1$)-simplex. A virtually equivalent criterion is given by Buck et al.~\cite{BIRS} (Corollary~1). The case $n=2$ has a remarkable consequence concerning positive braids.

\begin{theorem} The fibre surface of positive braid knots is obtained from the trivial open book by iterated positive trefoil plumbing.
\end{theorem}

At first sight, this seems to be a consequence of Corollary~3 in Giroux-Goodman~\cite{GN}. However, the discussion there only implies that positive braid knots are connected to the trivial knot by a sequence of positive trefoil plumbing and deplumbing operations. We conclude the introduction with an application of Theorem~2, the precise meaning of which we will explain in Section~5.

\begin{corollary} The fibre surface of a positive braid knot of genus $g$ can be untwisted by $g$ ribbon twists.
\end{corollary}

In the last section, we present a criterion by Hironaka~\cite{Hi} obstructing the existence of $n$-plumbing summands. As an application, we detect the largest $n$-plumbing summand for fibre surfaces of all 3-stranded torus links.

\section{Positive $n$-plumbing}

Let $(\Sigma, \phi)$ be an open book with a positive $n$-plumbing summand. The monodromy $\phi$ can be written as a composition
$$\phi=D_n D_{n-1} \ldots D_1 \bar{\phi},$$
where $D_1,D_2, \ldots,D_n$ are right-handed Dehn twists along the core curves of $n$ Hopf bands, in the order of plumbing, and the support of $\bar{\phi}$ is disjoint from the $n$ attached $1$-handles. Let $I \subset \Sigma$ be an essential relative arc in the outermost $1$-handle. The iterated images $I,\phi(I),\phi^2(I),\ldots \phi^n(I) \subset \Sigma$ are easy to determine: they are indeed homologically independent and pairwise disjoint, as shown in Figure~2, for $n=4$. This proves one implication of Theorem~1, without any further assumption on $\phi$.
\begin{figure}[ht]
\scalebox{1.0}{\raisebox{-0pt}{$\vcenter{\hbox{\epsffile{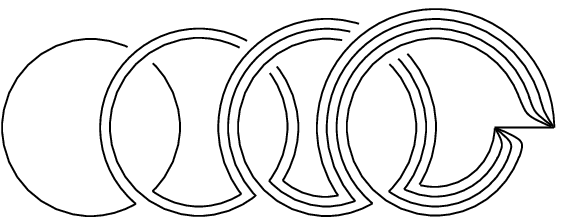}}}$}}
\caption{}
\end{figure}

For the reverse implication, we assume that $\phi$ is right-veering. Loosely speaking, this means that every relatively embedded interval $J \subset \Sigma$ is mapped to the right by $\phi$. Technically speaking, choose $\phi(J)$ to have minimal intersection number with $J$. Then the tangent vectors of $J$ and $\phi(J)$ are required to form a positive basis of $T \Sigma$, at their common boundary points (compare Figures~2 and~3). Right-veering diffeomorphisms were introduced by Honda et al.~\cite{HKM} in the context of contact geometry. One of their main results is that all closed tight contact $3$-manifolds are supported by right-veering open books. 

Let $I \subset \Sigma$ be a relatively embedded interval, such that $I, \phi(I)$, $\phi^2(I),\ldots,\phi^n(I)$ are homologically independent and have pairwise disjoint interiors. Setting
$$\gamma_k=\phi^{k-1}(I) \cup \phi^k(I),$$
we obtain $n$ homologically independent embedded circles $\gamma_1,\gamma_2,\ldots,\gamma_n \subset \Sigma$. After a suitable isotopy, these circles form a chain, i.e. the only intersection points arise from pairs of consecutive circles. Here the right-veering property is absolutely essential. This is illustrated in Figure~3 for $n=3$, where the bold and thin lines represent $I, \phi(I), \phi^2(I), \phi^3(I)$ and $\gamma_1,\gamma_2,\gamma_3$, respectively (the vertical line being $I$).
\begin{figure}[ht]
\scalebox{1.0}{\raisebox{-0pt}{$\vcenter{\hbox{\epsffile{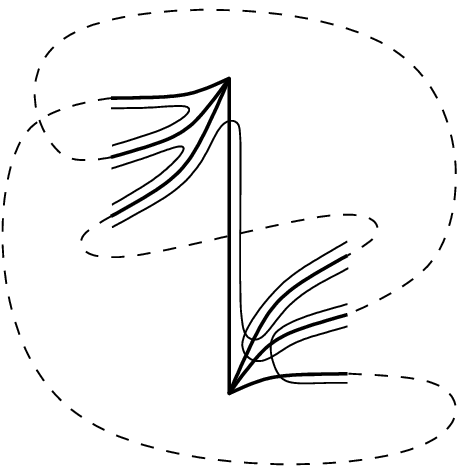}}}$}}
\caption{}
\end{figure}

Let $D_1,D_2,\ldots,D_n$ be right-handed Dehn twists along the curves $\gamma_1,\gamma_2,\ldots,\gamma_n$. Then the composition
$$\bar{\phi}=D_1^{-1} D_2^{-1} \ldots D_n^{-1} \phi$$
leaves the arcs $I,\phi(I),\phi^2(I),\ldots \phi^{n-1}(I)$ invariant. In particular, $\bar{\phi}$ is defined on the surface $\Sigma$ cut open along $I,\phi(I),\phi^2(I),\ldots \phi^{n-1}(I)$, which is again connected, due to the homological assumption. We conclude that $(\Sigma, \phi)$ contains a positive $n$-plumbing summand.

\section{The fibre and monodromy of positive braids}

A braid is positive, if it contains positive crossings only. Standard generators $\sigma_i$ are used to denote positive crossings. Braid words are read from top to bottom, in accordance with the second author's convention~\cite{De}. We use brick diagrams as an alternative notation for positive braids, as well as their fibre surfaces, see Figure~4. The diagram on the right shows the fibre surface of the braid $\sigma_3 \sigma_1 \sigma_2 \sigma_2 \sigma_3 \sigma_1 \sigma_2 \sigma_1$, an embedded orientable surface whose boundary is the closure of that braid and which naturally retracts on the brick diagram. It is instructive to verify that the isotopy type of the fibre surface is invariant under the braid relations $\sigma_{i} \sigma_{i+1} \sigma_{i} =\sigma_{i+1} \sigma_{i} \sigma_{i+1}$.
\begin{figure}[ht]
\scalebox{1.0}{\raisebox{-0pt}{$\vcenter{\hbox{\epsffile{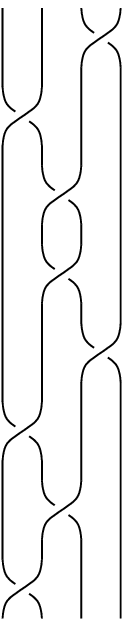}}}$}} \qquad \quad
\scalebox{1.0}{\raisebox{-0pt}{$\vcenter{\hbox{\epsffile{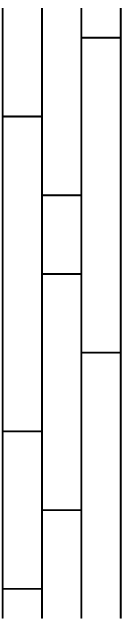}}}$}} \qquad \quad
\scalebox{1.0}{\raisebox{-0pt}{$\vcenter{\hbox{\epsffile{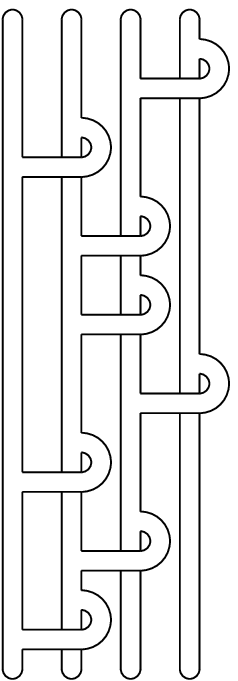}}}$}}
\caption{}
\end{figure}

Throughout this paper, we will only consider reduced braids, i.e. braids with at least two crossings between all strands. In other words, every column of the brick diagram contains at least one rectangle. Every rectangle defines a simple closed curve whose neighbourhood in the fibre surface is a positive Hopf band. The left- and topmost rectangle defines a Hopf band which can be deplumbed, since it lies on top of the surface (compare Figure~4).
On the level of braids, this amounts to removing one crossing on the top left. Therefore, fibre surfaces of reduced positive braids admit a positive Hopf plumbing structure. This well-known fact is due to Stallings~\cite{St}.

As explained in detail by Dehornoy~\cite{De}, the monodromy can be written as a product of right-handed Dehn twists, one for each rectangle, in the order of plumbing, i.e. from the bottom right to the top left, filling up columns from the right to the left. We will make use of this in the proof of Theorem~2.

\section{The plumbing structure of positive braids}

The proof of Theorem~2 relies on the following well-known fact, which can be found at various places in the literature, e.g. in Franks-Williams~\cite{FW}.

\begin{lemma} \quad
Every positive braid whose closure is a non-trivial link contains the square of a generator, $\sigma_m^2$, up to braid relations.
\end{lemma}

Here braid relations are understood in a generalised sense, including conjugation and Markov moves.

\begin{proof}[Proof of Lemma~1] Let $\beta$ be a positive braid representing a non-trivial link. By the braid relation, we may remove all instances of $\sigma_{1} \sigma_{2} \sigma_{1}$ in $\beta$. If there is only one single crossing of type $\sigma_1$ left, we remove it by a Markov move. Otherwise $\beta$ contains a section of the form $\sigma_1 w \sigma_1$ with either no $\sigma_2$ or at least two $\sigma_2$'s in $w$. In the first case we are done; in the second case $\beta$ contains a section of the form $\sigma_2 w \sigma_2$ with no $\sigma_1$ in $w$ and either no $\sigma_3$ or at least two $\sigma_3$'s in $w$. This argument eventually leads to a proof of the statement.
\end{proof}

The fibre surface associated with the braid $\sigma_1^2$ is a positive Hopf band. More generally, let $\beta$ be a positive braid that contains a square $\sigma_m^2$, then the corresponding fibre surface contains a Hopf band which can actually be deplumbed, since it lies on top of the surface. This provides an alternative proof that fibre surfaces of positive braids admit a positive Hopf plumbing structure. The monodromy of positive braid links is right-veering, since it can be written as a product of right-handed Dehn twists~\cite{HKM}. Therefore, we may apply the criterion of Theorem~1 to the open books associated with positive braids. Instead of looking at iterates of an arc $I \subset \Sigma$ under the monodromy~$\phi$, we will consider a simple closed curve~$R$ of the form $I \cup \phi(I)$ and verify that $\phi(R)$ does not intersect $I$. 

\begin{proof}[Proof of Theorem~2]
Let $\beta$ be a positive braid representing a non-trivial knot. By Lemma~1, we may assume that $\beta$ contains a square $\sigma_m^2$. After a suitable conjugation, we may further assume that $\beta$ starts with
$$\sigma_m^2 \sigma_{m-1}^{n_{m-1}} \sigma_{m-2}^{n_{m-2}} \ldots \sigma_1^{n_1},$$
for some numbers $n_1, \ldots, n_{m-1}>0$.
Moreover, $\beta$ contains at least one more generator $\sigma_m$, since it represents a knot.

Let $I \subset \Sigma$ be an essential relative arc in the uppermost band and let $R=I \cup \phi(I)$ be the simple closed curve running once through the Hopf band associated with $\sigma_m^2$. In order to determine the image $\phi(R)$, we decompose the monodromy as
$$\phi=\phi_1 \phi_2 \ldots \phi_l,$$
where $l$ is the number of columns of the brick diagram, and each $\phi_k$ is a product of right-handed Dehn twists in the $k$-th column, from the bottom to the top. The last $l-m$ factors $\phi_{m+1}, \ldots, \phi_l$ do not affect the curve~$R$. Indeed, a Dehn twist `along a rectangle' affects another rectangle, if and only if their curves intersect. This happens precisely when the two rectangles are arranged as in the braids
$\sigma_1^3$, $\sigma_1 \sigma_2 \sigma_1 \sigma_2$, or $\sigma_2 \sigma_1 \sigma_2 \sigma_1$. The effect of the monodromy on a single rectangle is described in detail in~\cite{De}.
The map $\phi_m$ shifts the curve~$R$ down, as shown on the left of Figure~5.
\begin{figure}[ht]
\scalebox{0.9}{\raisebox{-0pt}{$\vcenter{\hbox{\epsffile{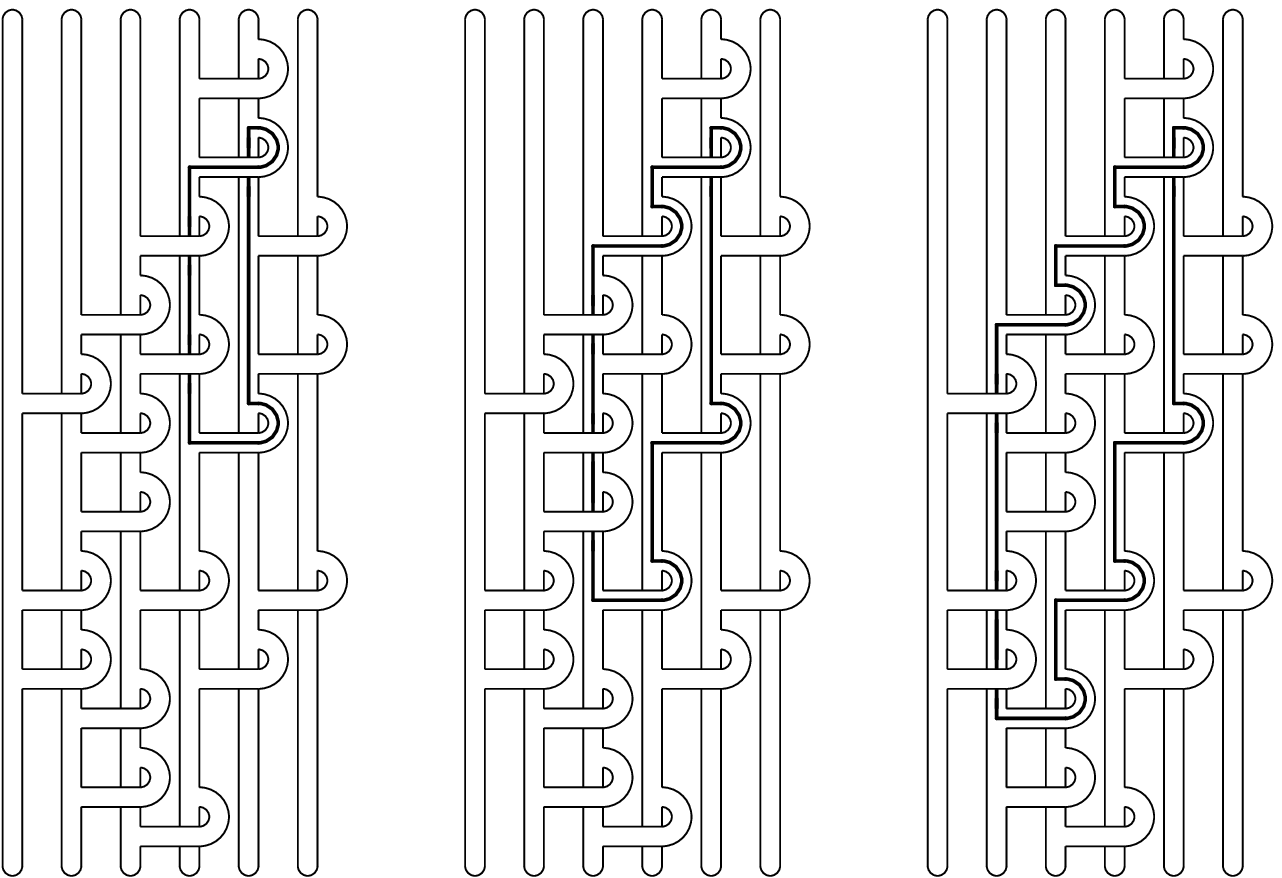}}}$}}
\caption{}
\end{figure}

The map $\phi_{m-1}$ adds a certain union of rectangles to its left. All subsequent factors $\phi_k$ ($k<m$) add one layer of rectangles to the curve $\phi_{k+1} \ldots \phi_l(R)$, as long as their union forms a staircase. Figure~5 illustrates this process for the braid
$$\sigma_4^2 \sigma_3 \sigma_2 \sigma_1 \sigma_5^2 \sigma_3 \sigma_4 \sigma_2^2 \sigma_5 \sigma_3 \sigma_4 \sigma_1^2 \sigma_2^2 \sigma_3.$$
The curves $\phi_4 \phi_5(R)$, $\phi_3 \phi_4 \phi_5(R)$ and $\phi_2 \phi_3 \phi_4 \phi_5(R)$ are shown from left to right.
In that example, the last factor $\phi_1$ does not have any effect on the curve $\phi_2 \phi_3 \phi_4 \phi_5(R)$.

We observe that $\phi(R)$ does not intersect the arc $I$. Now Theorem~1 allows us to apply trefoil deplumbing, which amounts to cutting $\Sigma$ along the arcs $I$ and $\phi(I)$, i.e. cutting the upper two bands associated with $\sigma_m^2$. On the level of braids, this simply means removing $\sigma_m^2$. It is instructive to verify that a neighbourhood of the theta graph $R \cup \phi(R) \subset \Sigma$ is indeed isotopic to the fibre surface of a positive trefoil knot (see again Figure~5). We conclude by induction. 
\end{proof}

\section{Untwisting positive braids}

A ribbon twist is an operation on Seifert surfaces that inserts a full twist into a ribbon, as shown in Figure~6. Two Seifert surfaces are twist equivalent if they are related by a finite number of ribbon twists. Generic Seifert surfaces are unlikely to be twist equivalent, since ribbon twists preserve the homotopy type of the surface complement in $\R^3$. Here we will show that fibre surfaces of positive braid knots can be untwisted in a natural way. A Seifert surface of genus $g$ is trivial if it admits $g$ disjoint non-isotopic compression discs on either side.\footnote{This makes sense even though a Seifert surface does not disconnect $\R^3$.} In particular, the complement of a trivial Seifert surface is homeomorphic to a handlebody of genus $2g$. A simple example for $g=1$ is depicted in Figure~6, on the right. 
\begin{figure}[ht]
\scalebox{0.6}{\raisebox{-0pt}{$\vcenter{\hbox{\epsffile{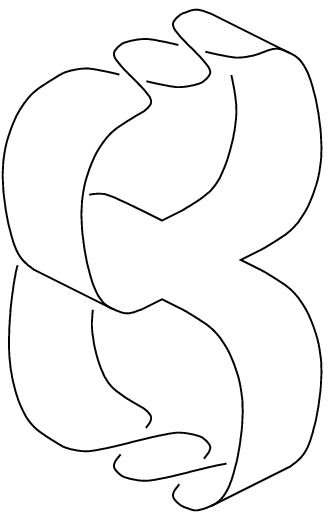}}}$}} $\quad \longrightarrow \quad$
\scalebox{0.6}{\raisebox{-0pt}{$\vcenter{\hbox{\epsffile{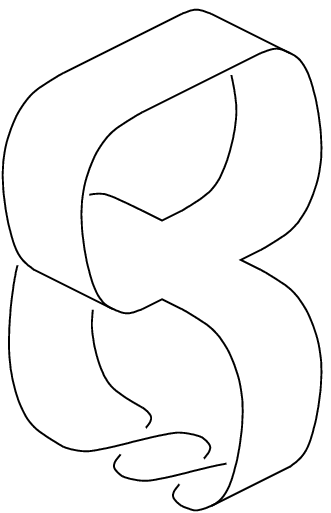}}}$}} $\qquad \cong \qquad$
\scalebox{0.6}{\raisebox{-0pt}{$\vcenter{\hbox{\epsffile{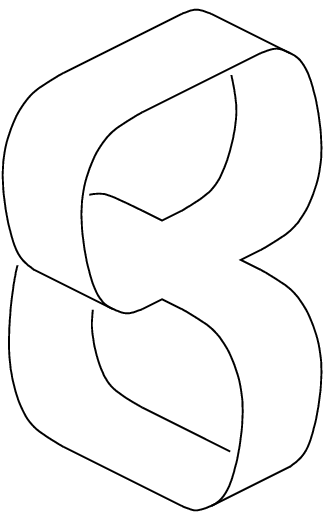}}}$}} 
\caption{}
\end{figure}

The sequence of diagrams in that figure shows that the fibre surface of the positive trefoil knot can be untwisted in one move. The last isotopy is best seen by dragging one end of the Hopf band once around the untwisted annulus.
In fact, the same sequence of diagrams, combined with the isotopy of Figure~7, demonstrates that positive trefoil plumbing can be expressed as a connected sum with a compressible torus, followed by a single ribbon twist. As a consequence, a Seifert surface $\Sigma \subset \R^3$ that is obtained by $g$-times iterated positive trefoil plumbing can be untwisted by $g$ ribbon twists. Moreover, the number $g$ is minimal, since the unknotting number of the boundary knot $\partial \Sigma \subset \R^3$ is not smaller than $g$ (see Rudolph~\cite{Ru}). This proves Corollary~1.
\begin{figure}[ht]
\scalebox{0.6}{\raisebox{-0pt}{$\vcenter{\hbox{\epsffile{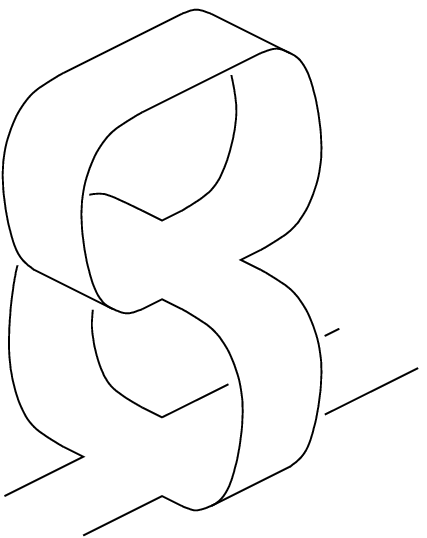}}}$}} $\qquad \cong \qquad$
\scalebox{0.6}{\raisebox{-0pt}{$\vcenter{\hbox{\epsffile{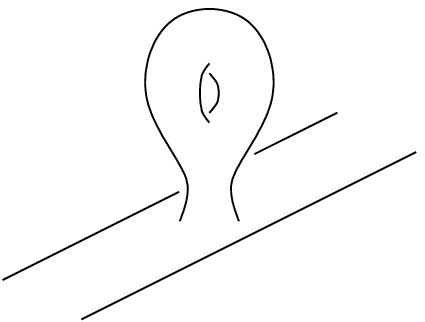}}}$}} 
\caption{}
\end{figure}

We conclude this section with a basic question: which fibre surfaces of genus $g$ in $\R^3$ can be untwisted by $g$ ribbon twists?

\section{Torus links and Hironaka's criterion}

The purpose of this section is to detect large $n$-plumbing summands in fibre surfaces of torus links. Using a result of Hironaka on the behaviour of the Alexander polynomial under iterated trefoil plumbing~\cite{Hi}, we will determine the largest $n$-plumbing summand for fibre surfaces of all 3-stranded torus knots. We define a torus link of type $T(p,q)$ as the closure of the braid
$$(\sigma_1 \sigma_2 \ldots \sigma_{p-1})^q.$$
The fibre surface of the torus knot $T(4,3)$ is shown in Figure~8. As explained in Section~3, the left- and topmost rectangle defines a Hopf band which can be deplumbed. Let $R$ be the core curve of that Hopf band. As in the proof of Theorem~2, we decompose the monodromy as $\phi=\phi_1 \phi_2 \ldots \phi_{p-1}$, where $\phi_k$ denotes the product of Dehn twists in the $k$-th column, from the bottom to the top. We claim that for all $k \leq p-2$, the curve $\phi^k(R)$ runs once around the uppermost rectangle in the $(k+1)$-th column. Indeed, let us suppose this is true for some $k<p-2$, and let us write $R_i$ ($1 \leq i \leq p-1$) for the curve defined by the uppermost rectangle in the $i$-th column (in particular, $R_{k+1}=\phi^k(R)$).
The first factor of the monodromy $\phi$ that affects the curve $R_{k+1}$ is $\phi_{k+2}$. In fact, $\phi_{k+2}(R_{k+1})$ is a curve running once around the union of $R_{k+1}$ and $R_{k+2}$. The next factor $\phi_{k+1}$ then removes $R_{k+1}$ from that curve, thus
$$\phi_{k+1} \phi_{k+2}(R_{k+1})=R_{k+2}.$$
The subsequent factors $\phi_1 \phi_2 \ldots \phi_k$ do not affect the curve $R_{k+2}$. Figure~8 illustrates this process for the torus knot $T(4,3)$ and $k=1$. For a more detailed analysis of the monodromy of Lorenz links, including torus links, we refer the reader to Section~2 of~\cite{De}.
\begin{figure}[ht]
\scalebox{1.0}{\raisebox{-0pt}{$\vcenter{\hbox{\epsffile{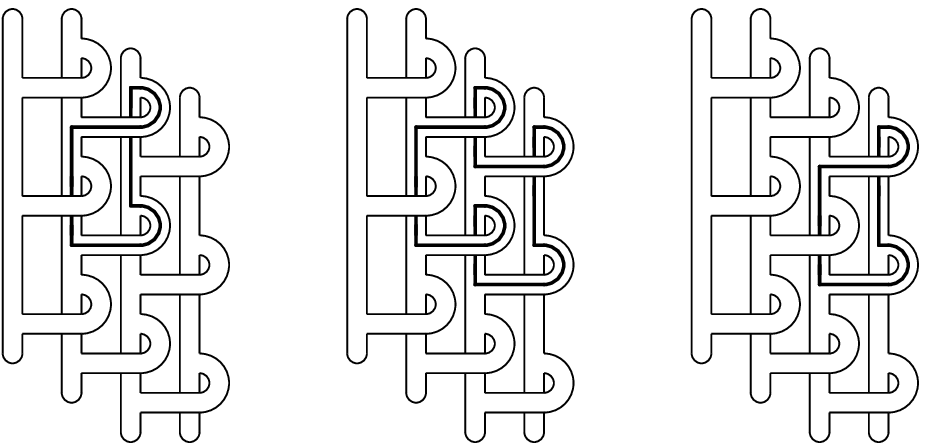}}}$}} 
\caption{}
\end{figure}

Keeping in mind that $R=I \cup \phi(I)$, where $I \subset \Sigma$ is an essential relative arc in the uppermost band, we deduce that the family of arcs $I,\phi(I),\phi^2(I),\ldots,\phi^{p-1}(I)$ is pairwise disjoint. As a consequence, the fibre surface $\Sigma(p,q)$ of the torus link $T(p,q)$ contains a $(p-1)$-plumbing summand. In the special case $q=2$, we obtain the entire surface as a  $(p-1)$-plumbing summand, as expected. In general, the fibre surface obtained by deplumbing a $(p-1)$-summand from $\Sigma(p,q)$ is not the fibre surface of a torus link (and not even the fibre surface of a positive braid link, for that matter).

In case a fibre surface $\Sigma$ admits an $n$-plumbing summand, we naturally obtain an embedded subsurface of type $\Sigma(2,n)$ in $\Sigma$. With this in mind, we will look for a large subsurface of type $\Sigma(2,n)$ in $\Sigma(3,m)$. Applying the braid relation three times, we obtain
$$(\sigma_1 \sigma_2)^6=\sigma_1^2 \sigma_2 \sigma_1^3 \sigma_2 \sigma_1^3 \sigma_2 \sigma_1.$$
In particular, the fibre surface $\Sigma(3,6)$ contains $\Sigma(2,9)$ as a subsurface (compare Figure~9). Roughly speaking, $\Sigma(3,m)$ contains a surface of type $\Sigma(2,n)$ which accounts for three quarters of the genus. However, this bound cannot be achieved by a plumbing summand.

\begin{proposition} Let $k$ be a natural number.
\begin{enumerate}
\item The fibre surface $\Sigma(3,3k+1)$ of the torus knot $T(3,3k+1)$ admits an $n$-plumbing summand, if and only if $n \leq 3k$. 
\item The fibre surface $\Sigma(3,3k+2)$ of the torus knot $T(3,3k+2)$ admits an $n$-plumbing summand, if and only if $n \leq 3k+2$. 
\end{enumerate}
\end{proposition}

\begin{figure}[ht]
\scalebox{1.0}{\raisebox{-0pt}{$\vcenter{\hbox{\epsffile{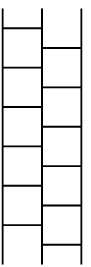}}}$}} $\qquad \cong \qquad$
\scalebox{1.0}{\raisebox{-0pt}{$\vcenter{\hbox{\epsffile{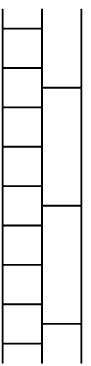}}}$}}
\caption{}
\end{figure}

By the above discussion, we are basically left to prove that $\Sigma(3,q)$ does not admit an $n$-plumbing summand for $n \geq q$ resp. $q+1$. This is a direct consequence of Hironaka's criterion, which involves the Alexander polynomial $\Delta_K(t)$ of links.

\begin{theorem}[\cite{Hi}, Theorem~9] Let $\Sigma_n$ be obtained by plumbing a $\Sigma(2,n)$-summand on the fibre surface $\Sigma$ of a link $K$ with $r$ components, and let $K_n$ be the boundary link of $\Sigma_n$. Then there exists a polynomial $P(t) \in \Z[t]$ of degree $d=\deg(\Delta_K(t))$, such that
$$\Delta_{K_n}(t)=\frac{t^n P(t)+(-1)^{r+n} t^d P(1/t)}{t+1}.$$
\end{theorem}

\begin{remark} Plumbing a $\Sigma(2,n)$-summand corresponds to positive $(n-1)$-plumbing in our convention, since $\Sigma(2,n)$ is a plumbing of $n-1$ Hopf bands.
\end{remark}

In order to appreciate this criterion, let us first analyse what it says about 2-stranded torus knots. The Alexander polynomial of torus knots is calculated in most textbooks about knot theory, e.g. in Rolfsen~\cite{Ro}:
$$\Delta_{T(p,q)}=\frac{(t^{pq}-1)(t-1)}{(t^p-1)(t^q-1)}.$$
The special case $p=2$ boils down to $\frac{t^q+1}{t+1}$. Comparing this with Hironaka's criterion, we obtain $P(t)=1$ and $n=q$, in accordance with the fact that $\Sigma(2,q)$ is itself a plumbing summand, and the remaining surface is a disc.

\begin{proof}[Proof of Proposition~1] \quad

\noindent
i) $q=3k+1$.
An elementary calculation shows
$$\Delta_{T(3,3k+1)}(t)=\sum_{i=0}^{2k} t^{3i}-\sum_{j=0}^{k-1}\left( t^{1+3j}+t^{6k-1-3j} \right).$$
Multiplying this with $(t+1)$ yields
$$1+t^{6k+1}+\sum_{j=0}^{k-1}\left( -t^{2+3j}+t^{3+3j}+t^{3k+1+3j}-t^{3k+2+3j} \right).$$
Comparing this again with Hironaka's criterion, we obtain
$$P(t)=\sum_{j=0}^{k-1}\left( t^{3j}-t^{1+3j} \right) +t^{3k},$$
thus $\deg(P)=3k$ and $n=3k+1$. For example,
\begin{align*}
(t+1) \Delta_{T(3,7)}(t) &= (t+1)(1-t+t^3-t^4+t^6-t^8+t^9-t^{11}+t^{12}) \\
&=1-t^2+t^3-t^5+t^6+t^7-t^8+t^{10}-t^{11}+t^{13},
\end{align*}
$P(t)=1-t+t^3-t^4+t^6$ and $n=7$.
We conclude that $\Sigma(2,3k+1)$ is the largest possible plumbing summand of the fibre surface $\Sigma(3,3k+1)$. In our convention, $\Sigma(3,3k+1)$ admits at most a positive $3k$-plumbing summand. Moreover, this bound is attained, by the discussion preceding Proposition~1.

\smallskip
\noindent
ii) $q=3k+2$. A similar calculation yields
$$\Delta_{T(3,3k+2)}(t)=-\sum_{i=0}^{2k} t^{1+3i}+\sum_{j=0}^{k}\left( t^{3j}+t^{6k+2-3j} \right),$$
$$P(t)=\sum_{j=0}^{k-1}\left( t^{3j}-t^{1+3j} \right) +t^{3k},$$
$\deg(P)=3k$ and $n=3k+3$. This leaves room for a $q$-plumbing summand, rather than $q-1$. In order to detect this summand, we work with the 3-stranded diagram of $T(3,3k+2)$. As usual, let $R=I \cup \phi(I)$ be the curve running once around the top left rectangle. A careful inspection of Figure~10 reveals that the third power of the monodromy $\phi=\phi_1 \phi_2$ shifts the curve $R$ down by three rectangles. Repeating this process $k$ times, we arrive at the bottom left rectangle. The latter is mapped to the bottom right rectangle under $\phi$. As a consequence, we obtain a family of pairwise disjoint arcs $I,\phi(I),\phi^2(I),\ldots,\phi^{3k+2}(I)$, in turn a $3k+2$-plumbing summand.
\begin{figure}[ht]
\scalebox{0.9}{\raisebox{-0pt}{$\vcenter{\hbox{\epsffile{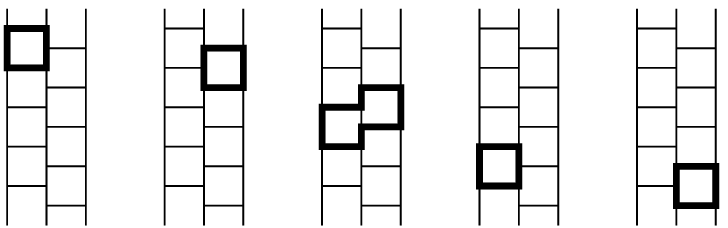}}}$}} 
\caption{}
\end{figure}
\end{proof}

The second item of Proposition~1 admits a generalisation to higher index torus knots.

\begin{proposition} Let $k$ and $p$ be natural numbers, $p$ odd. Then the fibre surface $\Sigma(p,pk+2)$ contains a $\Sigma(2,\frac{1}{2}p(p-1)k+p)$-plumbing summand.
\end{proposition}

\begin{remark} The genus of the fibre surface $\Sigma(p,q)$ of a torus knot $T(p,q)$ is $\frac{1}{2}(p-1)(q-1)$. Applying this to the knots of Proposition~2, we see that the largest plumbing summand of $\Sigma(p,pk+2)$ fills up more than half of the genus (the precise ratio is $\frac{1}{2}+\frac{1}{2(pk+1)}$ and is again given by Hironaka's criterion).
\end{remark}

\begin{figure}[ht]
\scalebox{0.9}{\raisebox{-0pt}{$\vcenter{\hbox{\epsffile{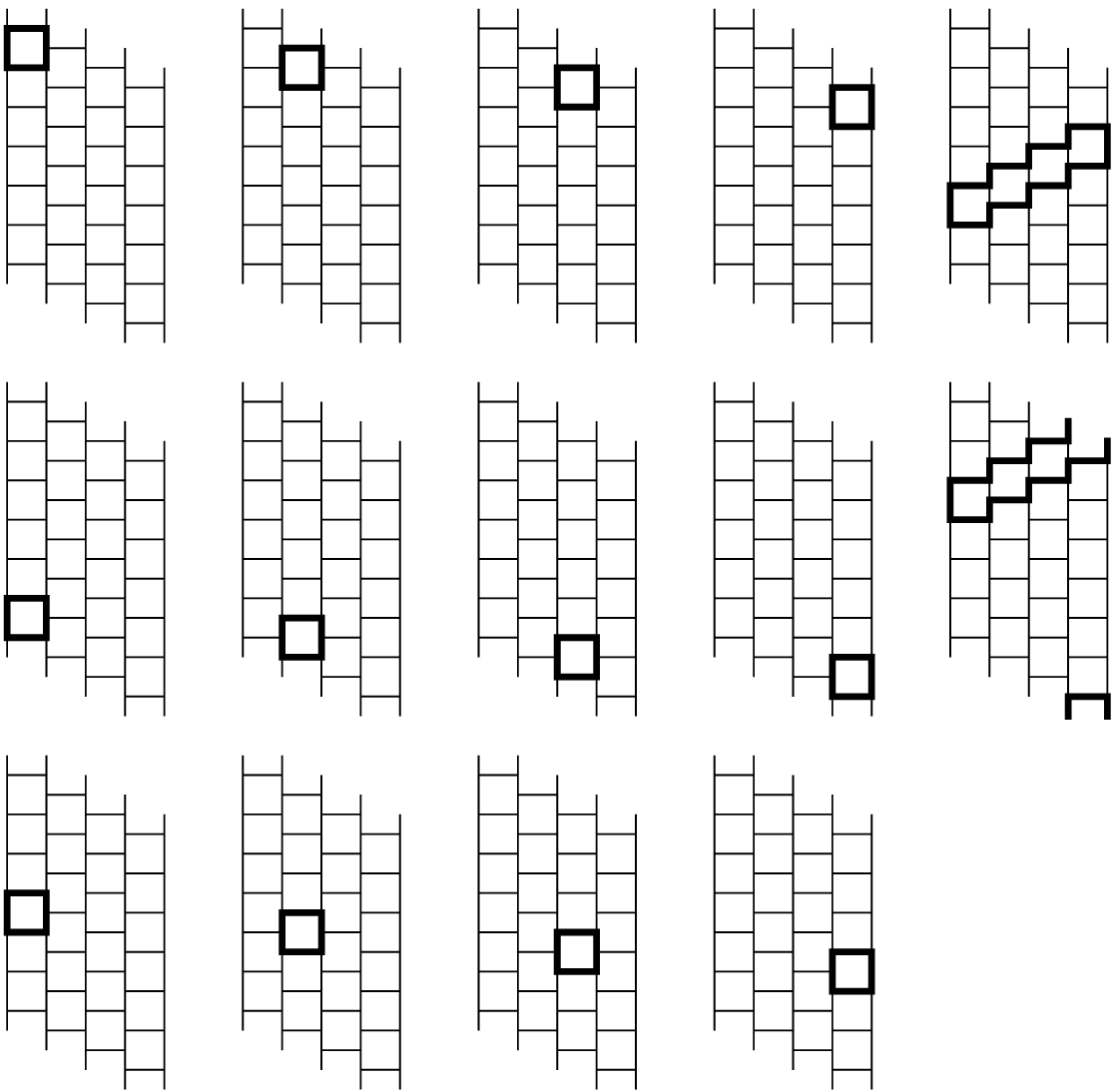}}}$}} 
\caption{}
\end{figure}

\begin{proof}[Proof of Proposition~2]
As in the above proofs, it is all about tracing the image of the top left rectangle $R=I \cup \phi(I)$ under iterates of the monodromy $\phi$. The $p$-th power of the monodromy shifts $R$ down by $p$ rectangles. This process can be repeated $\frac{1}{2}(p-1)k$ times, without ever intersecting the initial curve $I$. Another $p-2$ iterations map the last curve $\phi^{\frac{1}{2}p(p-1)k}(R)$ to a rectangle in the right column. The whole process is shown in Figure~11, for the torus knot $T(5,7)$. Summing up, the first $\frac{1}{2}p(p-1)k+p-1$ iterates of the arc $I$  are pairwise disjoint. This gives rise to a plumbing summand of type $\Sigma(2,\frac{1}{2}p(p-1)k+p)$, as required. 
\end{proof}

In general, Hironaka's bound allows much larger $n$-plumbing summands than we were able to detect.

\bigskip
\noindent
Universit\"at Bern, Sidlerstrasse 5, CH-3012 Bern, Switzerland

\bigskip
\noindent
\texttt{sebastian.baader@math.unibe.ch}

\smallskip
\noindent
\texttt{pierre.dehornoy@math.unibe.ch}

\end{document}